\documentclass[10pt]{article}
\usepackage{amssymb,amsbsy,amsmath,amsfonts,amssymb,amscd}
\usepackage{latexsym,euscript,exscale,epic,eepic,epsfig}
\usepackage[english,francais]{babel}

\newtheorem{e-proposition}[theorem]{Proposition}

\newtheorem{e-definition}[theorem]{Definition\rm}

\newtheorem{theoreme}{Th\'eor\`eme}%[section]

\newtheorem{proposition}[theoreme]{Proposition}

\def\adots{\mathinner{\mkern2mu\raise1pt\hbox{.}
\mkern3mu\raise4pt\hbox{.}\mkern1mu\raise7pt\hbox{.}}}
\def\<{\langle\,}
\def\>{\,\rangle}

\def\shuff#1#2{\mathbin{
      \hbox{\vbox{
        \hbox{\vrule
              \hskip#2
              \vrule height#1 width 0pt
               }%
        \hrule}%
             \vbox{
        \hbox{\vrule
              \hskip#2
              \vrule height#1 width 0pt
               \vrule }%
        \hrule}%
}}}
\def\shuffl{{\mathchoice{\shuff{7pt}{3.5pt}}%
                        {\shuff{6pt}{3pt}}%
                        {\shuff{4pt}{2pt}}%
                        {\shuff{3pt}{1.5pt}}}}%
\def\shuffle{\, \shuffl \,}

\def\Sylv{{\rm Sylv}}
\def\SG{{\mathfrak S}}
\def\Sh{{\rm Sh}}
\def\Sym{{\bf Sym}}
\def\FQSym{{\bf FQSym}}
\def\FSym{{\bf FSym}}
\def\PBT{{\bf PBT}}
\def\F{{\bf F}}
\def\G{{\bf G}}
\def\PP{{\bf P}}
\def\QQ{{\bf Q}}
\def\pp{{\cal P}}
\def\qq{{\cal Q}}
\def\SS{{\bf S}}
\def\T{{\cal T}}

\def\C{{\mathbb C}}

\def\std{{\rm std\, }}

\def\K{{\mathbb K}}

\def\maj{{\rm maj\,}}
\def\today{\number\day \space\ifcase\month\or
 Janvier\or F\'evrier\or Mars\or Avril\or Mai\or Juin\or
 Juillet\or Ao\^ut\or Septembre\or Octobre\or Novembre\or D\'ecembre\fi
 \space\number\year}

\begin{document}
\selectlanguage{francais}
\title{%
Un analogue du mono\"\i de plaxique\\
  pour les arbres binaires de recherche 
}
\author{%
Florent HIVERT,\ 
Jean-Christophe NOVELLI,\ 
Jean-Yves THIBON
}
\maketitle
\thispagestyle{empty}
%%%%%%%%%%%%%%%%%%%%%%%%%%%%%%%%%%%%%%%%%%%%%%%%%%%%%%%%%%%%
%%%  R\'esum\'e  %%%
%%%%%%%%%%%%%%%%%%%%
\begin{abstract}{%

Nous introduisons une structure de mono\"\i de sur un ensemble
d'arbres binaires \'etiquet\'es, par un proc\'ed\'e analogue
\`a la construction du mono\"\i de plaxique. Nous en d\'eduisons
une nouvelle approche de l'alg\`ebre des arbres binaires
de Loday-Ronco.
}\end{abstract}
%%%%%%%%%%%%%%%%%%%%%%%%%%%%%%%%%%%%%%%%%%%%%%%%%%%%%%%%%%%%
%%%  Le titre en Anglais  %%%
%%%%%%%%%%%%%%%%%%%%%%%%%%%%%
\selectlanguage{english}
\begin{center}{%
{\bf An analogue of the plactic monoid for binary search trees}
}\end{center}
%%%%%%%%%%%%%%%%%%%%%%%%%%%%%%%%%%%%%%%%%%%%%%%%%%%%%%%%%%%%
%%%  Abstract  %%%
%%%%%%%%%%%%%%%%%%
\begin{abstract}{%
We introduce a monoid structure on a certain set of labelled
binary trees, by a process similar to the construction of the
plactic monoid. This leads to a new interpretation of the algebra of planar
binary trees of Loday-Ronco.
}\end{abstract}

%%%%%%%%%%%%%%%%%%%%%%%%%%%%%%%%%%%%%%%%%%%%%%%%%%%%%%%%%%%%
%%%  Abridged English version  %%%
%%%%%%%%%%%%%%%%%%%%%%%%%%%%%%%%%%
%\AEv

%La version anglaise abr\'eg\'ee.

\par\medskip\centerline{\rule{2cm}{0.2mm}}\medskip
\setcounter{section}{0}
\selectlanguage{francais}
%%%%%%%%%%%%%%%%%%%%%%%%%%%%%%%%%%%%%%%%%%%%%%%%%%%%%%%%%%%%
%%%  Texte principal (en Francais)  %%%
%%%%%%%%%%%%%%%%%%%%%%%%%%%%%%%%%%%%%%%
%Le texte principal (en Fran\c{c}ais).

Il existe certaines analogies entre la combinatoire des arbres binaires et
celle des tableaux de Young. Par exemple, les \'etiquetages croissants
d'un arbre binaire donn\'e et les tableaux de Young standard d'une forme
donn\'ee sont tous deux d\'enombr\'es par une  \og formule des \'equerres\fg
 \cite{FRT,St,Kn} dont on conna\^\i t dans les deux cas des $q$-analogues naturels
\cite{St,BW1}.
On sait aussi qu'il est possible de construire, au moyen de la correspondance
de Robinson-Schensted, une alg\`ebre de Hopf ayant pour base les tableaux
de Young standard  \cite{PR}, et 
qu'une r\'ealisation naturelle de cette alg\`ebre au moyen de polyn\^omes
non commutatifs, 
implique une preuve combinatoire \'eclairante (en une ligne) de la r\`egle de
Littlewood-Richardson \cite{Loth2,DHT,NCSF6}.
Dans cette r\'ealisation,
chaque tableau $t$ de forme $\lambda$ est un polyn\^ome
homog\`ene de degr\'e $n=|\lambda|$ dont l'image commutative est la fonction de Schur
$s_\lambda$.

R\'ecemment, Loday et Ronco ont introduit une alg\`ebre de Hopf
ayant pour base les arbres binaires planaires \cite{LR1,LR2}.
Cette alg\`ebre peut, comme la pr\'ec\'edente, 
se r\'ealiser dans une alg\`ebre associative libre \cite{DHT,NCSF6},
chaque arbre binaire complet \`a $n$ sommets internes \'etant
repr\'esent\'e par un polyn\^ome homog\`ene de degr\'e $n$
en variables non commutatives.
En r\'ealit\'e, ces deux alg\`ebres sont, par d\'efinition, des sous-big\`ebres
de l'alg\`ebre de Hopf des permutations de Malvenuto-Reutenauer \cite{MvR},
laquelle a \'et\'e r\'ealis\'ee dans \cite{DHT,NCSF6} comme
l'alg\`ebre des fonctions quasi-sym\'etriques libres $\FQSym$ dont nous
rappelons ci-dessous la d\'efinition.
Notons au passage que toutes les alg\`ebres en question sont libres \cite{PR,LR1,MvR},
mais que cette propri\'et\'e n'intervient pas dans la construction de leurs
r\'ealisations polynomiales.
 
Ce sont ces r\'ealisations qui permettent de faire  appara\^\i tre l'alg\`ebre
des tableaux (ou fonctions sym\'etriques libres) et l'alg\`ebre
des arbres binaires comme deux cas particuliers d'une m\^eme construction,
reposant sur l'existence d'une correspondance de type Robinson-Schensted
et d'un mono\"\i de de type plaxique. 
Cette construction nous appara\^\i t d'autant plus fondamentale qu'un
troisi\`eme exemple 
entre dans ce cadre, celui du couple d'alg\`ebres de Hopf en
dualit\'e $(\Sym,QSym)$ (fonctions sym\'etriques non commutatives
et fonctions quasi-sym\'etriques), lequel correspond au mono\"\i de hypoplaxique
\cite{NCSF4}.

Nous commencerons par rappeler quelques d\'efinitions.
Nos notations seront celles de \cite{Loth2} et de \cite{DHT,NCSF6}.
Soit $A=\{a_1<a_2<\cdots\}$ un alphabet d\'enombrable totalement ordonn\'e.
Les polyn\^omes non commutatifs
$\F_\sigma(A)=\sum_{\std(w)=\sigma^{-1}}w$, o\`u $\sigma$ est une permutation et
o\`u $w$ parcourt les mots de standardis\'e $\sigma^{-1}$ forment une base d'une
sous alg\`ebre de $\K\<A\>$ not\'ee $\FQSym_\K$ ($\K$ etant un anneau commutatif
quelconque) et simplement $\FQSym$ si  $\K=\C$ \cite{DHT}.
On pose \'egalement $\G_\sigma=\F_{\sigma^{-1}}$ et on d\'efinit un
produit scalaire par $\<\F_\sigma,\G_\tau\>=\delta_{\sigma\tau}$.
 
Soit $w\mapsto(P(w),Q(w))$ la correspondance de Robinson-Schensted 
usuelle ({\it cf.} \cite{Loth2}).
Pour un tableau standard $t$ de forme $\lambda$, on d\'efinit la fonction de
Schur libre $\SS_t$ comme la somme  $\SS_t=\sum_{P(\sigma)=t}\F_\sigma$. Ainsi qu'il
est montr\'e dans \cite{Loth2,NCSF6}, les $\SS_t$ forment une base d'une sous-alg\`ebre
$\FSym$ de $\FQSym$, et cet \'enonc\'e peut \^etre vu comme un raffinement
de la r\`egle de Littlewood-Richardson, qu'il implique imm\'ediatement.
De plus, $\FSym$ est une sous-big\`ebre de $\FQSym$. 
 
Nous allons maintenant donner une construction de l'alg\`ebre $\PBT$ des arbres
binaires planaires de Loday-Ronco \cite{LR1} enti\`erement analogue \`a
celle de $\FSym$. Commen\c cons par l'analogue appropri\'e de la correspondance
de Robinson-Schensted. Elle se d\'efinit \`a partir de l'arbre binaire
d\'ecroissant (ou arbre tournoi)
$\T(\sigma)$ associ\'e \`a une permutation $\sigma\in\SG_n$.
C'est un arbre binaire (incomplet) \`a $n$ sommets, num\'erot\'es de $1$ \`a $n$.
La racine est num\'erot\'ee $n$ (le maximum), et si, en tant que mot, $\sigma=unv$,
alors
le sous-arbre gauche est $\T(u)$ et le sous-arbre droit $\T(v)$ ($u$ et $v$ \'etant
vus comme des permutations de leurs r\'earrangements croissants).                   
On notera $[\T(\sigma)]$ la forme de cet arbre (on oublie les \'etiquettes).
 
Pour un mot $w\in A^*$, on posera $\qq(w)=\T((\std w)^{-1})$, et $\pp(w)$ sera l'arbre
obtenu en rempla\c cant chaque \'etiquette $i$ de $\qq(w)$ par la $i$-\`eme
lettre de $w$. Par exemple, si $w=bacaabca$, $\std w=\sigma=51723684$,
$\sigma^{-1}=24581637$ et
$$
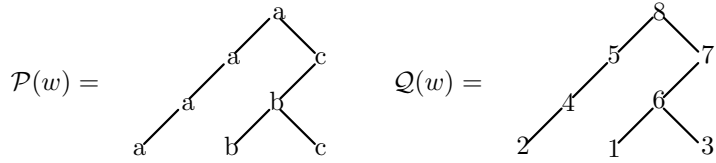
\pp(w) =
\vcenter{\hbox{
\setlength{\unitlength}{0.0002in}
\begingroup\makeatletter\ifx\SetFigFont\undefined%
\gdef\SetFigFont#1#2#3#4#5{%
  \reset@font\fontsize{#1}{#2pt}%
  \fontfamily{#3}\fontseries{#4}\fontshape{#5}%
  \selectfont}%
\fi\endgroup%
{\renewcommand{\dashlinestretch}{30}
\begin{picture}(4999,3825)(0,-10)
\put(4800,0){\makebox(0,0)[lb] c}
\thicklines
\path(3600,3600)(2700,2700)
\path(2400,2400)(1500,1500)
\path(1200,1200)(300,300)
\path(4800,2400)(3900,1500)
\path(3600,1200)(2700,300)
\path(3900,1200)(4800,300)
\put(3675,3600){\makebox(0,0)[lb] a}
\put(2475,2400){\makebox(0,0)[lb] a}
\put(1275,1200){\makebox(0,0)[lb] a}
\put(0,0){\makebox(0,0)[lb] a}
\put(2400,0){\makebox(0,0)[lb] b}
\put(3600,1200){\makebox(0,0)[lb] b}
\put(4800,2400){\makebox(0,0)[lb] c}
\path(3900,3600)(4800,2700)
\end{picture}
}
}
}
\qquad
\qq(w)=
\vcenter{\hbox{
\setlength{\unitlength}{0.0002in}
\begingroup\makeatletter\ifx\SetFigFont\undefined%
\gdef\SetFigFont#1#2#3#4#5{%
  \reset@font\fontsize{#1}{#2pt}%
  \fontfamily{#3}\fontseries{#4}\fontshape{#5}%
  \selectfont}%
\fi\endgroup%
{\renewcommand{\dashlinestretch}{30}
\begin{picture}(5090,3915)(0,-10)
\put(4875,75){\makebox(0,0)[lb] 3}
\thicklines
\path(3600,3675)(2700,2775)
\path(2400,2475)(1500,1575)
\path(1200,1275)(300,375)
\path(3600,1275)(2700,375)
\path(4800,2475)(3900,1575)
\path(3899,1276)(4799,376)
\put(3600,3600){\makebox(0,0)[lb] 8}
\put(2400,2475){\makebox(0,0)[lb] 5}
\put(1200,1200){\makebox(0,0)[lb] 4}
\put(0,75){\makebox(0,0)[lb] 2}
\put(2400,0){\makebox(0,0)[lb] 1}
\put(3600,1275){\makebox(0,0)[lb] 6}
\put(4875,2475){\makebox(0,0)[lb] 7}
\path(3900,3675)(4800,2775)
\end{picture}
}
}
}
$$
Ces arbres peuvent se calculer au moyen d'un algorithme d'insertion \`a la Schensted.
Pour construire $\pp(w)$, on lit les lettres de $w$ \`a partir de la droite. La
derni\`ere est plac\'ee \`a la racine, puis chaque lettre est 
ins\'er\'ee r\'ecursivement  dans le
sous-arbre gauche ou droit selon qu'elle est $\le$ ou $>$ \`a
la racine. L'arbre $\qq(w)$ indique l'ordre (inverse) de cr\'eation
des sommets. Ainsi, l'\'etiquetage de $\pp(w)$ est croissant au sens large
dans le sens $\nearrow$ et au sens strict dans le sens $\searrow$, et le parcours
infixe  de l'arbre produit le r\'earrangement croissant de $w$. Il s'agit l\`a
d'un algorithme de tri classique \cite{Kn}, et on reconna\^\i t dans $\pp(w)$ un
arbre binaire de recherche. Quant \`a $\qq(w)$, c'est un arbre tournoi.
Dans le cas o\`u $w=\sigma$ est une permutation, la seule information contenue
dans $\pp(\sigma)$ est sa forme, et on l'identifie \`a l'arbre non \'etiquet\'e
$[\T(\sigma^{-1})]$.
 
\begin{theoreme}\label{thsylv}
La relation d'\'equivalence $\sim$ sur le mono\"\i de libre $A^*$ d\'efinie
par
$$
u\sim v\ \Longleftrightarrow \ \pp(u)=\pp(v)
$$
co\"\i ncide avec la congruence engendr\'ee par les relations
\begin{equation}
zxuy  \equiv xzuy\,,\ x\le y< z\in A\,,\ u\in A^*\,.
\end{equation}
\end{theoreme}
Cette congruence sera appel\'ee {\it congruence sylvestre} et le mono\"\i de
quotient $\Sylv(A)=A^*/\equiv$, le {\it mono\"\i de sylvestre}.
 
La congruence sylvestre est aux arbres ce que la congruence plaxique
(engendr\'ee par les relations de Knuth) est aux tableaux \cite{LS,Loth2}. 
Dans le cas du mono\"\i de plaxique, les classes ont des repr\'esentants canoniques
(les lectures par lignes des tableaux). Il en est de m\^eme ici :
 
\begin{proposition}
Soit $T$ un arbre binaire de recherche, et soit $w_T$ le mot obtenu en
effectuant son parcours postfixe (parcours en profondeur  commen\c cant
par la gauche, o\`u l'on \'ecrit la lettre \'etiquetant un noeud quand on le visite
pour la derni\`ere fois). Alors, $\pp(w_T)=T$, et $w_T$ est 
minimal pour l'ordre lexicographique dans sa classe sylvestre.
De plus, en it\'erant \`a partir de $w_T$
les r\`egles de r\'e\'ecriture $xzuy\rightarrow zxuy$, obtenues
en orientant les relations sylvestres, on engendre toute la classe de $w_T$.
\end{proposition}

Les mots de la forme $w_T$ seront appel\'es {\it mots-arbres} (ou plus simplement
{\it arbres}).
                                                                                    
Pour un arbre binaire $T$ non \'etiquet\'e (identifi\'e \`a un $\pp$-symbole de
permutation), introduisons le polyn\^ome
$$
\PP_T(A)=\sum_{\pp(\sigma)=T}\F_\sigma\,.
$$
En utilisant le fait que la congruence sylvestre est compatible \`a la restriction
aux intervalles de l'alphabet, et en raisonnant comme dans 
\cite{NCSF6}, prop. 3.12,on en d\'eduit l'identit\'e suivante.
 
\begin{theoreme}
Soient $T'$ et $T''$ deux arbres binaires non \'etiquet\'es. Alors,
$$
\PP_{T'}\PP_{T''}=\sum_{T\in \Sh(T',T'')}\PP_T\,,
$$
o\`u $\Sh(T',T'')$ d\'esigne l'ensemble des arbres $T$ tels que
$w_T$ appara\^\i t dans le
produit de m\'elange $u\shuffle v$ o\`u $u=w_{T'}$ et o\`u $v=w_{T''}[k]$ est le mot-arbre
de $T''$ d\'ecal\'e du nombre $k$ de sommets de $T'$.
\end{theoreme}
 
Cet \'enonc\'e est enti\`erement analogue \`a la
\og r\`egle de Littlewood-Richardson libre\fg
  pr\'esent\'ee dans \cite{Loth2,NCSF6}. Il fournit une nouvelle construction
de l'alg\`ebre des arbres binaires de Loday-Ronco. Dans \cite{LR2}, le
produit $\PP_{T'}\PP_{T''}$ est d\'ecrit au moyen d'un ordre sur les arbres
binaires (notons au passage que Loday et Ronco utilisent des arbres binaires
complets, dont nous ne conservons ici que les noeuds internes).
Cet ordre peut s'obtenir \`a partir de l'ordre faible du groupe sym\'etrique :
c'est la restriction de celui-ci aux permutations qui sont des mots-arbres.
Remarquons que les r\'esultats de Bj\"orner et Wachs \cite{BW2} entra\^\i nent que les
classes sont des intervalles de l'ordre faible.
 
Ainsi, l'ensemble $\Sh(T',T'')$ est un intervalle de l'ordre de Loday-Ronco.
Dans le cas des tableaux, il est possible de d\'efinir un ordre similaire,
quotient de l'ordre faible par les relations de Knuth,
dont les ensembles $\Sh(t',t'')$ de \cite{Loth2,NCSF6}  sont 
\'egalement des intervalles.
 
Le cardinal d'une classe plaxique de permutations est \'egal au nombre
de tableaux standard d'une certaine forme, donn\'e par la c\'el\`ebre
formule des \'equerres \cite{FRT}, dont on conna\^\i t un $q$-analogue
\cite{St} qui prend en compte l'indice majeur des permutations.
De la m\^eme mani\`ere, le d\'enombrement par indice majeur de la classe
sylvestre correpondant \`a l'arbre $T$ \`a $n$ sommets est donn\'e
par la sp\'ecialisation $(q)_n\PP_T(1,q,q^2,\cdots)$, \'egale, d'apr\`es \cite{BW1}
\`a
\begin{equation}
\sum_{\pp(\sigma)=T}q^{\maj(\sigma)}={[n]_q! \over \prod_{\circ\in T}
q^{\delta_\circ}[h_\circ]_q}\,,
\end{equation}
o\`u, pour un sommet $\circ$ de $T$, $h_\circ$ est le nombre de sommets
du sous-arbre dont il est racine, et $\delta_\circ$ celui de son
sous-arbre droit.

La congruence sylvestre permet \'egalement une description
simple du dual de Hopf de l'alg\`ebre des arbres binaires :

\begin{theoreme}
Le dual $\PBT^*$ de $\PBT$ est isomorphe \`a l'image
de $\FQSym$ par la projection canonique
$$
\pi :\ \C\<A\>\longrightarrow \C[\Sylv(A)]\simeq \C\<A\>/\equiv\,.
$$
La base duale $\QQ_T$ de $\PP_T$ est $\QQ_T=\pi(\G_{w_T})$,
o\`u $w_T$ est la permutation canoniquement associ\'ee \`a
l'arbre $T$.
\end{theoreme}
Le dual de $\FSym$ admet une description similaire \cite{PR,NCSF7}.
Ces r\'esultats sont des cons\'equences de la \og formule de Cauchy libre\fg
$$
\sum_{\std(u)=\std(v)^{-1}}u\otimes v
=\sum_\sigma {\bf U}_\sigma\otimes {\bf V_\sigma}
$$
pour tout couple $({\bf U},{\bf V})$ de bases adjointes
de $\FQSym$.

Pour calculer les \'el\'ements primitifs de $\PBT$ et de $\PBT^*$
(ou de $\FSym$ et de $\FSym^*$), on peut proc\'eder comme
dans \cite{NCSF6}. On part d'une base multiplicative d'un c\^ot\'e.
Les \'el\'ements de la base duale correspondant aux ind\'ecomposables
de notre base multiplicative forment alors une base d'\'el\'ements
primitifs du dual. Le calcul explicite se ram\`ene \`a celui de la
fonction de M\"obius de l'ordre appropri\'e.
Dans la cas de $\PBT$, et plus g\'en\'eralement des alg\`ebres dendriformes
libres, les \'el\'ements primitifs ont \'et\'e d\'etermin\'es dans \cite{Ro}.
Toutefois, le proc\'ed\'e propos\'e ci-dessus conduit \`a une caract\'erisation
diff\'erente.

Finalement, remarquons qu'on peut construire un couple de graphes
gradu\'es $(\Gamma,\Gamma^*)$ en dualit\'e au sens de Fomin \cite{Fo},
dont les sommets de degr\'e $n$ sont les arbres binaires
\`a $n$ sommets. Dans $\Gamma$, on a une ar\^ete entre $T$ et
$T'$ si $\PP_{T'}$ appara\^\i t dans le d\'eveloppement
de $\PP_T\PP_\bullet$. Dand $\Gamma^*$, cette m\^eme ar\^ete
sera pr\'esente si $\QQ_{T'}$ appara\^\i t dans $\QQ_T\QQ_\bullet$.
La correspondance sylvestre est alors la correspondance de Fomin associ\'ee
\`a ce couple de graphes.

\medskip
%%%%%%%%%%%%%%%%%%%%%%%%%%%%%%%%%%%%%%%%%%%%%%%%%%%%%%%%%%%%
%%%  Remerciements  %%%
%%%%%%%%%%%%%%%%%%%%%%%
{\footnotesize Nous tenons \`a remercier Jean-Louis Loday et Mar\'\i a Ronco pour leurs commentaires.}
%%%%%%%%%%%%%%%%%%%%%%%%%%%%%%%%%%%%%%%%%%%%%%%%%%%%%%%%%%%%
%%%  Bibliographie %%%
%%%%%%%%%%%%%%%%%%%%%%

%
\bigskip

\sc
Institut Gaspard Monge,\\
Universit\'e de Marne-la-Vall\'ee,\\
77454 Marne-la-Vall\'ee cedex,\\
FRANCE\\

\footnotesize
\begin{thebibliography}{99}
\selectlanguage{english}
%
\bibitem{BW1} A. Bj\"orner, M. Wachs, $q$-Hook length formulas for forests,
J. Combin. Theory Ser. A 52 (1989), 165--187.
%
\bibitem{BW2} A. Bj\"orner, M. Wachs, Permutation statistics and linear extensions
of posets,
J. Combin. Theory Ser. A 58 (1991), 85--114.
%
\bibitem{DHT} G. Duchamp, F. Hivert, J.-Y. Thibon,
Une g\'en\'eralisation des fonctions quasi-sym\'etriques et des fonctions
sym\'etriques non commutatives,
C. R. Acad. Sci. Paris S\'er. I Math. 328 (1999), no. 12, 1113--1116.
%
%
\bibitem{NCSF6}  G. Duchamp, F. Hivert, J.-Y. Thibon,
Noncommutative symmetric functions VI: free quasi-symmetric functions and related algebras,
Internat. J. Alg. Comput. (\`a para\^\i tre).
%
\bibitem{NCSF7} G. Duchamp, F. Hivert, J.-Y. Thibon, 
Noncommutative symmetric functions VII (en pr\'eparation).
%
\bibitem{Fo} S. Fomin, Duality of graded graphs, J. Alg. Combinatorics
3 (1994), 357--404.
%
\bibitem{FRT} J.S. Frame, G. de B. Robinson, R.M. Thrall,
The hook graphs of the symmetric groups, Canadian J. Math. 6, (1954). 316--324.
%
\bibitem{Kn} D.E. Knuth, The art of computer programming, vol.3: Soting
and searching, Addison-Wesley, 1973.
%
\bibitem{NCSF4} D. Krob, J.-Y. Thibon, Noncommutative symmetric functions IV:
Quantum linear groups and Hecke algebras at $q=0$,
J. Algebraic Combin. 6 (1997), no. 4, 339--376.
%
\bibitem{LS} A. Lascoux, M.-P. Sch\"utzenberger,
Le monoïde plaxique, Noncommutative structures in algebra and geometric combinatorics 
(Naples, 1978), pp. 129--156, Quad. Ricerca Sci., 109, CNR,
Rome, 1981.
%
\bibitem{LR1} J.-L. Loday, M.O. Ronco,
Hopf algebra of the planar binary trees, Adv. Math. 139 (1998), no. 2, 293--309.
%
\bibitem{LR2} J.-L. Loday, M.O. Ronco,
Order structure on the algebra of permutations and of planar binary trees,
J.  Algebraic Combin. 15  (2002), no. 3,  253-270.
%
\bibitem{Loth2} M. Lothaire, Algebraic combinatorics on words, Cambridge
University press, Cambridge, 2002.
%
\bibitem{MvR} C. Malvenuto, C. Reutenauer, Duality between quasi-symmetric functions and Solomon descent
algebra, J. Algebra 177 (1995), 967--892.
%
\bibitem{PR} S. Poirier, C. Reutenauer,
Alg\`ebres de Hopf de tableaux,
Ann. Sci. Math. Qu\'ebec 19 (1995), no. 1, 79--90.
%
\bibitem{Ro} M.O. Ronco, Primitive elements in a free dendriform Hopf algebra, 
Contemp. Maths. Vol. 267 (2000), 245--264. 
%
\bibitem{St} R.P. Stanley,
Ordered structures and partitions, 
Memoirs of the American Mathematical Society, No. 119, 
American Mathematical Society, Providence, R.I., 1972.
%                                 
%%
\end{thebibliography}
\end{document}